\documentclass[12pt, a4paper]{article}

\usepackage{amssymb}
\usepackage{amsmath}
\usepackage{amscd}
\usepackage{theorem}
\usepackage{mathrsfs}

\usepackage{graphicx}
\usepackage[T1]{fontenc}
\usepackage{ae}
\usepackage{pictexwd,dcpic}

\theoremstyle{change}

\newtheorem{thm}{Theorem.}
\newtheorem{lem}[thm]{Lemma.}

\newcommand{\N}{\mathbb{N}}

\newcommand{\C}{\mathbb{C}}

\renewcommand{\phi}{\varphi}

\newcommand{\proof}{\textbf{Proof. }}
\newcommand{\proofend}{\hfill $\Box$}

\begin{document}

\begin{center}
{\textbf{\large Toeplitz operators with pluriharmonic symbol}}
\end{center}

\vspace{1cm}

\centerline{J\"org Eschmeier and Sebastian Langend\"orfer}

\vspace{.8cm}

\begin{center}
\parbox{12cm}{\small Let $m \geq 1$ be an integer and let $H_m(\mathbb B)$ be 
the analytic functional Hilbert space on the unit ball $\mathbb B \subset \mathbb C^n$ given by the
reproducing kernel $K_m(z,w) = (1 - \langle z,w \rangle)^{-m}$. We prove that Toeplitz 
operators with pluriharmonic symbol on $H_m(\mathbb B)$ can be characterized by an
algebraic identity which extends the classical Brown-Halmos characterization of
Toeplitz operators on the Hardy space of the unit disc as well as corresponding results of
Louhichi and Olofsson for Toeplitz operators with harmonic symbol on weighted Bergman spaces
of the disc.
 
\vspace{0.5cm}

}

\end{center}
\vspace{1cm}

\centerline{\textbf{\S1 \, Introduction}}

\vspace{.5cm}

A result of Brown and Halmos \cite{BH} from 1963 shows that an operator $T \in L(H^2(\mathbb T))$
on the Hardy space of the unit disc is a Toeplitz operator $T_f = P_{H^2(\mathbb T)} M_f|_{H^2(\mathbb T)}$
with $L^{\infty}$-symbol $f \in L^{\infty}(\mathbb T)$ if and only if the operator $T$ satisfies
the algebraic identity $M_z^*TM_z = T$. In 2008 Louhichi and Olofsson \cite{LO} proved that an operator
$T \in L(A_m(\mathbb D))$ on the standard weighted Bergman space
\[
A_m(\mathbb D) = \{ f \in \mathcal O(\mathbb D); \| f \|^2 = \frac{m-1}{\pi} 
\int_{\mathbb D} |f(z)|^2 (1 - |z|^2)^{m-2} dz < \infty \}
\]
is a Toeplitz operator with harmonic symbol $f$ on $\mathbb D$ if and only if $T$ satisfies the algebraic
identity
\[
M_z^{\prime *}TM_z^{\prime} = \sum_{k=0}^{m-1} (-1)^k \binom{m}{k+1} M^k_z T M_z^{*k}.
\]
Here $M_z^{\prime} = M_z(M^*_zM_z)^{-1} \in L(A_m(\mathbb D))$ is the Cauchy dual of the multiplication
operator $M_z: \, A_m(\mathbb D) \rightarrow A_m(\mathbb D), g \mapsto zg$. 

The weighted Bergman space
$A_m(\mathbb D)$ is the analytic functional Hilbert space with reproducing kernel
$K_m(z,w) = (1 - z\overline{w})^{-m}$. The aim of the present note is to extend the above result from
\cite{LO} to the analytic functional Hilbert spaces $H_m(\mathbb B)$ on the unit ball $\mathbb B \subset \mathbb C^n$
given by the reproducing kernels $K_m(z,w) = (1 - \langle z,w \rangle)^{-m}$, where $m \geq 1$ is an
arbitrary positive integer. We show that an operator $T \in L(H_m(\mathbb B))$ is a Toeplitz operator with
pluriharmonic symbol $f$ on $\mathbb B$ if and only if the operator $T$ satisfies the algebaic identity
\[
M^{\prime *}_zTM^\prime_z = 
P_{{\rm Im} M^*_z}\Big(\oplus \sum^{m-1}_{k=0}(-1)^k \binom{m}{k+1} \sigma^k_{M_z}(T)\Big)P_{{\rm Im}M^*_z}.
\]
Here $M^*_z$ is the adjoint of the row operator $M_z: \, H_m(\mathbb B)^n \rightarrow H_m(\mathbb B), 
(g_i) \mapsto \sum_{1 \leq i \leq n} z_i g_i$, $M^\prime_z$ is a suitably defined Cauchy dual of
the multiplication tuple $M_z$ and the operator $\sigma_{M_z} \in L(H_m(\mathbb B))$ acts as
$\sigma_{M_z}(T) = \sum_{1 \leq i \leq n} M_{z_i}TM_{z_i}^*$ . 

For $m \in \{1, \ldots, n-1 \}$, the multiplication tuple $M_z$ is not subnormal and it is not immediately
obvious how to define Toeplitz operators with non-analytic symbols in these cases (see e.g. \cite{AK, Tch}). 
Using the fact that 
pluriharmonic functions on $\mathbb B$ admit a decomposition $f = g + \overline{h}$ with analytic
functions $g,h \in \mathcal O(\mathbb B)$, we suggest a natural definition of Toeplitz operators
with pluriharmonic symbol on analytic functional Hilbert spaces over $\mathbb B$.

For $m = 1$, the space $H_1(\mathbb B)$ is the Drury-Arveson space
on the unit ball and the above algebraic identity characterizing Toeplitz operators with pluriharmonic
symbol simplifies on $H_1(\mathbb B)$ to the identity
\[
M^*_z T M_z = P_{{\rm Im} M^*_z}(\oplus T)P_{{\rm Im}M^*_z}.
\]
If in addition $n = 1$, then $M_z^*$ is onto and the identity reduces to the classical 
Brown-Halmos condition. For $m = n$, the space
$H_n(\mathbb B)$ is the Hardy space 
\[
H_n(\mathbb B) = \{ f \in \mathcal O(\mathbb B),\ \| f \|^2 = 
\sup\limits_{0<r< 1} \int\limits_{\partial \mathbb B}| f(r \xi)|^2 d\sigma(\xi)<\infty \},
\]
while for $m \geq n+1$, the space $H_m(\mathbb B)$ is the weighted Bergman space
\[
H_m(\mathbb B) =
\{ f \in \mathcal O(\mathbb B); \; \| f \|^2 = 
\int_{\mathbb B} |f|^2 d\mu_m < \infty \},
\]
of all analytic functions that are square integrable with respect to the measure
$\mu_m = \frac{(m-1)!}{(m-n-1)! \pi^n} (1 - |z|^2)^{m-n-1} dz$, where $dz$ is the 
Lebesgue measure on $\mathbb C^n$. In these cases our definition of
Toeplitz operators with pluriharmonic symbol agrees with the usual definition and our
characterization of Toeplitz operators with pluriharmonic symbol extends the 
one-dimensional results of Louhichi and Olofsson from \cite{LO}.

\vspace{1cm}

\centerline{\textbf{\S2 \, Results}}

\vspace{.5cm}

Let $H_m(\mathbb B)$ be the functional Hilbert space with reproducing kernel
\[
K_m:\mathbb B \times \mathbb B \rightarrow \mathbb C,\ K_m(z,w)=\frac{1}{(1-\langle z,w\rangle)^m}
= \sum_{k=0}^{\infty} \binom{m+k-1}{k} \langle z,w \rangle^k,
\]
where $m > 0$ is a positive integer. Then
\[
H_m(\mathbb B)=\lbrace f=\sum_{\alpha\in \mathbb N^n}f_\alpha z^\alpha \in \mathcal O(\mathbb B);\ \| f\|^2=\sum_{\alpha \in \mathbb N^n}\frac{|f_\alpha|^2}{\rho_m(\alpha)}< \infty\rbrace
\]
with $\rho_m(\alpha)=\frac{(m+|\alpha|-1)!}{\alpha!(m-1)!}$. 
The tuple
\[
M_z = (M_{z_1},\ldots,M_{z_n}) \in L(H_m(\mathbb B))^n
\]
consisting of the multiplication operators
$
M_{z_i}:H_m(\mathbb B)\rightarrow H_m(\mathbb B),\ f\mapsto z_i f,
$
with the coordinate functions is well defined and its Koszul complex
\[
K^{\cdot}(M_z,H_m(\mathbb B))\stackrel{\epsilon_\lambda}{\longrightarrow}\mathbb C \longrightarrow 0
\]
augmented by the point evaluation $\epsilon_\lambda : H_m(\mathbb B)\rightarrow \mathbb C,\ f \mapsto f(\lambda)$, is exact. 
In particular, the last map in the Koszul complex
\[
H_m(\mathbb B)^n\stackrel{M_z}{\longrightarrow} H_m(\mathbb B),\ (f_i)^n_{i=1} \mapsto \sum^n_{i=1} z_i f_i
\]
has closed range $M_z H_m(\mathbb B)^n = \lbrace f \in H_m(\mathbb B);\ f(0)=0\rbrace$.
The above properties of the functional Hilbert spaces $H_m(\mathbb B)$ are well known and follow for
instance from the results in Section 4 of \cite{Guo}.

In the following we write $M_z : \,  H_m(\mathbb B)^n \rightarrow H_m(\mathbb B)$ for the row multiplication
defined as above, and we write $M^*_z : H_m(\mathbb B) \rightarrow H_m(\mathbb B)^n,\ f\mapsto (M^*_{z_i}f)^n_{i=1}$, 
for its adjoint. Since the row operator $M_z$ has closed range, the
operator $M^*_z M_z: {\rm Im}\ M^*_z \rightarrow {\rm Im} M^*_z$ is invertible. We use the notation 
$(M^*_z M_z)^{-1}$ for its inverse. The space $H_m(\mathbb B)$ admits the 
orthogonal decomposition
\[
H_m(\mathbb B) = \bigcirc\!\!\!\!\!\!\!\perp^\infty_{k=0}\mathbb H_k
\]
into the spaces $\mathbb H_k \subset \mathbb C[z]$ consisting of all homogeneous polynomials of degree $k$. 
For each function $f = \sum_{\alpha \in \mathbb N^n} f_{\alpha} z^\alpha \in H_m(\mathbb B)$, 
its homogeneous expansion
\[
f = \sum^{\infty}_{k=0}f_k \; \mbox{  with  } \; f_k = \sum_{|\alpha|=k}f_\alpha z^\alpha,
\]
coincides with the decomposition of $f$ as an element of the orthogonal direct sum $H_m(\mathbb B)=\bigcirc\!\!\!\!\!\!\!\perp^\infty_{k=0} \mathbb H_k$. Let us denote by $\delta,\Delta:H_m(\mathbb B)
\rightarrow H_m(\mathbb B)$ the invertible diagonal operators acting as
\[
\delta (\sum^\infty_{k=0}f_k ) = f_0 + \sum^\infty_{k=1}\frac{m+k-1}{k}f_k
\]
and $\Delta\left(\sum^\infty_{k=0}f_k\right) = \sum^\infty_{k=0}\frac{m+k}{k+1}f_k$. One can show that
\[
M^*_z \delta f = (M^*_z M_z)^{-1}(M^*_z f)
\]
for all $f \in H_m(\mathbb B)$ and that the row operator $\delta M_z: H_m(\mathbb B)^n \rightarrow H_m(\mathbb B)$ 
defines a continuous linear extension of the operator
\[
M_z (M^*_z M_z)^{-1}: {\rm Im} M^*_z \rightarrow H_m(\mathbb B).
\]
(Lemma 1 in \cite{E}). The diagonal operators $\delta$ and $\Delta$ satisfy the intertwining relation 
$\delta M_z = M_z(\oplus \Delta)$ and $\Delta$ admits the representation
\[
\Delta = \sum^{m-1}_{j=0}(-1)^j \binom{m}{j+1} 
\sum_{|\alpha| = j} \gamma_{\alpha} M^{\alpha}_z M^{* \alpha}_z = 
\sum^{m-1}_{j=0}(-1)^j\binom{m}{j+1}\sigma^j_{M_z}(1_{H_m(\mathbb B)})
\]
(Lemma 3 in \cite{E}). Here $\gamma_{\alpha} = |\alpha|!/\alpha!$ for $\alpha \in \mathbb N^n$ and $\sigma_{M_z}\in L(H_m(\mathbb B))$ is defined by
\[
\sigma_{M_z}(X)=\sum^n_{i=1}M_{z_i}XM^\ast_{z_i}.
\]
In dimension $n=1$, the operator 
\[
M^\prime_z = \delta M_z : \, H_m(\mathbb B)^n \rightarrow H_m(\mathbb B)
\]
is usually called the Cauchy dual of $M_z$ (see e.g. \cite{Sh}).

We begin by showing that the powers of $M_z$ and $M^*_z$ satisfy a Brown-Halmos type condition.

\begin{lem}\label{powers}
For each multiindex $\gamma \in \mathbb N^n$, the identity
\[
M^{\prime *}_z (M^{\gamma}_z) M^\prime_z = 
P_{{\rm Im}M^*_z}\Big(\oplus\sum^{m-1}_{j=0}(-1)^j\binom{m}{j+1}\sigma^j_{M_z}(M^{\gamma}_z)\Big)P_{{\rm Im}M^*_z}
\]
holds.
\end{lem}

\proof
Note that ${\rm Ker} M_z = ({\rm Im}M^*_z)^\perp$. Hence we obtain\\
\begin{align*}
M^{\prime *}_z M^{\gamma}_z M^{\prime}_z 
&= P_{{\rm Im}M^*_z} \big(M^*_z\delta M^{\gamma}_z \delta M_z\big)P_{{\rm Im}M^*_z}\\
&=P_{{\rm Im}M^*_z} \big((M^*_z M_z)^{-1}M^*_z M^{\gamma}_z M_z(\oplus \Delta)\big)P_{{\rm Im}M^\ast_z}\\
&=P_{{\rm Im}M^*_z} \big((M^*_z M_z)^{-1}(M^*_z M_z)(\oplus M^{\gamma}_z \Delta \big) P_{{\rm Im}M^*_z}\\
&=P_{{\rm Im}M^*_z} \big(\oplus \sum^{m-1}_{j=0}(-1)^j\binom{m}{j+1}\sigma^j_{M_z}(M^{\gamma}_z) \big)P_{{\rm Im}M^*_z}.
\end{align*} 
Here we have used the identity $(M^\ast_zM_z)^{-1}(M^\ast_zM_z)=P_{{\rm Im}M^\ast_z}$. 
\proofend

By passing to adjoints we find that also the powers of $M^*_z$ satisfy the identity
\[
M^{\prime *}_z (M^{* \gamma}_z) M^{\prime}_z = 
P_{{\rm Im} M^*_z} \Big(\oplus\sum^{m-1}_{j=0}(-1)^j\binom{m}{j+1}\sigma^j_{M_z}(M^{* \gamma}_z)\Big)P_{{\rm Im}M^*_z}.
\]

Let us denote by $\mathcal T_{{\rm BH}} \subset L(H_m(\mathbb B))$ the set of all 
bounded linear operators on $H_m(\mathbb B)$ that satisfy the Brown-Halmos type condition
\[
M^{\prime \ast}_zTM^\prime_z = P_{{\rm Im}M^\ast_z}\Big(\oplus \sum^{m-1}_{j=0}(-1)^j\binom{m}{j+1}\sigma^j_{M_z}(T)\Big)P_{{\rm Im}M^\ast_z}.
\]
Our aim is to show that $\mathcal T_{{\rm BH}}$ consists precisely of all Toeplitz operators with pluriharmonic
symbol. To show that every operator $T$ in $\mathcal T_{{\rm BH}}$ is a Toeplitz operator we decompose $T$
into its homogeneous components.

For this purpose, we denote by $U: \mathbb R \rightarrow L(H_m(\mathbb B))$ the strongly continuous unitary 
operator group acting as $(U(t)f)(z) = f(e^{{\rm it}}z)$. Then $H_m(\mathbb B)$ admits the orthogonal decomposition
\[
H_m(\mathbb B) = \bigcirc\!\!\!\!\!\!\!\perp_{k\in \mathbb Z} H_k,
\]
where the spaces $H_k$ are the images of the orthogonal projections $P_k \in L(H_m(\mathbb B))$ defined by
\[
P_k = \frac{1}{2 \pi} \int^{2 \pi}_0 e^{-ikt} U(t)dt.
\]
Here the integrand is regarded as a continuous function with values in the locally convex space
$L(H_m(\mathbb B))$ equipped with the strong operator topology and the integral is defined as a
weak integral (see Theorem 3.27 in \cite{Ru2} and  Section 20.6(3) in \cite {K}). All operator-valued integrals
used in the following should be understood in this sense.

An application of Cauchy's integral theorem yields that $P_k = 0$ for $k < 0$.
For $t \in \mathbb R$ and $k\in \mathbb Z$, the identity $U(t)P_k = e^{ikt}P_k$ holds. The space $H_k$ 
consists precisely of all functions $f \in H_m(\mathbb B)$ with
\[
f(e^{{\rm it}}z) = e^{ikt}f(z)
\]
for all $t \in \mathbb R$ and all $z \in \mathbb B$. Thus $H_k = \mathbb H_k$ for $k \geq 0$.
For a bounded operator $T \in L(H_m(\mathbb B))$ and $k \in \mathbb Z$, we define its $k$th homogeneous 
component $T_k \in L(H_m(\mathbb B))$ by
\[
T_k = \frac{1}{2\pi}\int \limits^{2\pi}_0 e^{-ikt}U(t)TU(t)^* dt.
\]
Let $k \in \mathbb Z$. The $k$th homogeneous component of the adjoint $T^*$ of $T$ is given by
$(T^*)_k = (T_{-k})^*$. To check this identity it suffices to observe that
\begin{align*}
\langle (T^*)_k f,g \rangle &= \frac{1}{2\pi}\int \limits^{2\pi}_0 \langle e^{-ikt}U(t) T^* U(t)^* f,g \rangle dt\\
&= \frac{1}{2\pi}\int \limits^{2\pi}_0 \langle f, e^{ikt}U(t)TU(t)^* g \rangle dt = \langle f, T_{-k}g \rangle 
\end{align*}
for all $f,g \in H_m(\mathbb B)$.
Since $U(t)^*|H_j = e^{-ijt}1_{H_j}$, we obtain that
\[
T_kf = \frac{1}{2 \pi}\int\limits^{2\pi}_0 e^{-i(k+j)t}U(t)Tf\ dt = P_{k+j}Tf \in H_{j+k}
\]
for $f \in H_j$. Thus $T_k H_j \subset H_{j+k}$ for all $j \in \mathbb Z$. Any bounded operator 
on $H_m(\mathbb B)$ with this property will be called homogeneous of degree $k$.

A standard argument using the Fej\'er kernel $K_N(t) = \sum_{|k| \leq N} (1 - \frac{|k|}{N+1}) e^{ikt}$
as summability kernel (Lemma I.2.2 in \cite{Katz}) shows that, for each $f \in H_m(\mathbb B)$,
\[
Tf = \lim_{N \rightarrow \infty} \frac{1}{2\pi}\int \limits^{2\pi}_0 K_N(t) U(t)TU(t)^* f dt
   = \lim_{N \rightarrow \infty} \sum_{|k| \leq N} (1 - \frac{|k|}{N+1}) T_k f
\]

\begin{lem}  \label{inherit}
Let  $T \in \mathcal T_{{\rm BH}}$ be given. Then $T_k \in \mathcal T_{{\rm BH}}$ for all $k \in \mathbb Z$.
\end{lem}

\proof
A straightforward calculation shows that
\[
U(t)^* \delta M_z = e^{-it}\delta M_z(\oplus U(t)^*)
\]
and hence $M^*_z \delta U(t) = e^{it}(\oplus U(t))M^*_z \delta$ for 
$t \in \mathbb R$. But then 
\[
M^*_z \delta T_k \delta M_z = \frac{1}{2\pi}\int\limits^{2 \pi}_0 
e^{-ikt}(\oplus U(t))(M^*_z\delta T\delta M_z)(\oplus U(t)^*)dt
\]
for $k\in \mathbb Z$. Since $U(t)M_{z_{j}} = e^{it}M_{z_{j}}U(t)$ for $t \in \mathbb R$ and $j = 1,\ldots, n$, 
the space ${\rm Im}M^*_z$ is reducing for $\oplus U(t)$ and
\[
U(t)M^\alpha_zTM^{* \alpha}_z U(t)^* = M^\alpha_z U(t)TU(t)^* M^{*\alpha}_z \; \; (t \in \mathbb R,\alpha \in \mathbb N^n).
\]
Using the hypothesis that $T \in \mathcal T_{{\rm BH}}$, we find that $M^*_z \delta T_k \delta M_z$
is given by
\begin{align*}
& \quad \frac{1}{2\pi}\int\limits^{2\pi}_0 
e^{-ikt}(\oplus U(t))P_{{\rm Im}M^*_z}\Big(\oplus \sum^{m-1}_{j=0}
(-1)^j\binom{m}{j+1}\sigma^j_{M_z}(T)\Big)P_{{\rm Im}M^*_z}) (\oplus U(t)^*) dt\\
&= P_{{\rm Im}M^*_z}\Big(\frac{1}{2\pi}\int\limits^{2\pi}_0 
e^{-ikt}(\oplus \sum^{m-1}_{j=0} (-1)^j \binom{m}{j+1} \sigma^j_{M_z}\big(U(t)TU(t)^*)\big)dt\Big) P_{{\rm Im}M^*_z}\\
&= P_{{\rm Im}M^*_z}\Big(\oplus \sum\limits^{m-1}_{j=0}(-1)^j \binom{m}{j+1} \sigma^j_{M_z}(T_k)\Big) P_{{\rm Im}M^*_z}.
\end{align*}
Thus we have shown that $T_k\in \mathcal T_{{\rm BH}}$ for every $k\in \mathbb Z$.
\proofend

All operators $T \in \mathcal T_{{\rm BH}}$ that are homogeneous of non-negative degree act as mutliplication operators.

\begin{thm} \label{multi}
Let $T \in \mathcal T_{{\rm BH}}$ be homogeneous of degree $r \in \mathbb N$. 
Then $T$ acts as the multiplication operator
\[
Tf = (T1)f  \quad (f \in H_m(\mathbb B)).
\]
\end{thm}

\proof
Define $q = T1 \in \mathbb{H}_r$. 
We write $M_q: H_m(\mathbb B) \rightarrow H_m(\mathbb B), f \mapsto qf$, for 
the operator of multiplication with $q$
and show by induction on $k$ that $T = M_q$ on $\mathbb H_k$ for all $k \in \N$. \\
For $k = 0$, this is obvious. Suppose that the assertion has been proved 
for $j = 0,...,k$ and fix a polynomial $p \in \mathbb{H}_{k+1}$. By Lemma 1 
in \cite{E} we have
\begin{align*}
M_z^* \delta T \delta M_z (M_z^* p) &= M_z^* \delta T M_z (M_z^* M_z)^{-1} M_z^* p 
= M_z^* \delta T P_{\text{Im}M_z} p \\
&= M_z^* \delta T p 
= \frac{m+k+r}{k+r+1} M_z^* (Tp).
\end{align*}
Here we ave used that $p \in \mathbb{H}_{k+1} \subseteq \C^{\bot} = \text{Im}M_z$. \\
Using the induction hypothesis and Lemma 3 from \cite{E}, we find that
\begin{align*}
& P_{{\rm Im}M_z^*}\Big(\oplus \sum_{j=0}^{m-1} (-1)^j \binom{m}{j+1}\sigma_{M_z}^j(T) \Big)P_{{\rm Im}M_z^*} (M_z^* p) \\
= & P_{{\rm Im}M_z^*} \Big(\oplus \sum_{j=0}^{m-1} (-1)^j \binom{m}{j+1} \sum_{| \alpha | = j} 
\gamma_{\alpha} M_z^{\alpha}T M_z^{*\alpha}\Big) (M_{z}^* p)\\
= & P_{{\rm Im } M_z^*} \Big(\oplus M_q \sum_{j=0}^{m-1} (-1)^j 
\binom{m}{j+1}\sigma_{M_z}^j (1_{H_m(\mathbb{B})})\Big)(M_z^* p) \\
=& P_{{\rm Im}M_z^*} \big(\oplus M_q \Delta \big)(M_z^* p) 
=  \frac{m+k}{k+1} P_{{\rm Im}M_z^*} \big(\oplus M_q\big) (M_z^* p).
\end{align*}

By hypothesis
\[
\frac{m+k+r}{k+r+1} M_z^* Tp = \frac{m+k}{k+1} P_{{\rm Im}M_z^*} (\oplus M_q) M_z^* p.
\]

By applying the operator $M_z(M_z^* M_z)^{-1} = \delta M_z|_{{\rm Im}\, M^*_z}$ to both sides of this equation, and by
using the identities
\[
(M^*_z M_z)^{-1} (M^*_z M_z) = P_{{\rm Im}M_z^*}, \; M_z(M_z^* M_z)^{-1} M^*_z = P_{{\rm Im}M_z}, 
\]
we find that
\begin{align*}
\frac{m+k+r}{k+r+1} Tp &= M_z(M_z^* M_z)^{-1} (\frac{m+k+r}{k+r+1} M^*_z Tp)\\ 
&= M_z(M_z^* M_z)^{-1} (\frac{m+k}{k+1}  P_{{\rm Im}M_z^*} (\oplus M_q) M_z^* p)\\ 
&= \frac{m+k}{k+1} \delta M_z(M_z^* M_z)^{-1} (M_z^*M_z) (\oplus M_q) M_z^* p \\
&= \frac{m+k}{k+1} \delta M_z (\oplus M_q) M_z^* p\\
&= \frac{m+k}{k+1} \; \frac{m+k+r}{k+r+1} M_q M_z M^*_z p.
\end{align*}
Since $M_z M^*_z = \sum_{j=1}^{\infty} \frac{j}{m+j-1}P_{\mathbb H_j}$ (Proposition 4.3 in \cite{Guo}), we conclude that
\[
Tp = \frac{m+k}{k+1} M_q \frac{k+1}{m+k}p = M_q p.
\]
This observation completes the induction and hence the whole proof.
\proofend

In general it is not obviuos how to define Toeplitz 
operators with non-analytic symbols on the full range of all analytic 
Besov-Sobolev spaces $H_m(\mathbb B)$ $(m \geq 1)$. However, for the case of 
pluriharmonic symbols, this problem can easily be
overcome. To begin with, let us fix a function $f \in H_m(\mathbb B)$. Then 
the set $D_f=\lbrace u \in H_m(\mathbb B);\ fu\in H_m(\mathbb B)\rbrace \subset H_m(\mathbb B)$ is 
a dense linear subspace which contains $\mathbb C[z]$, and
\[
T_f:D_f\rightarrow H_m(\mathbb B),\ u \mapsto fu
\]
is a densely defined closed linear operator. For $g \in H_m(\mathbb B)$, let us denote by $g_\alpha=g^\alpha(0)/\alpha!$ its Taylor coefficients at $z=0$. Then for $\alpha \in \mathbb N^n$ and
$u \in D_f$,
\[
\langle fu,z^\alpha\rangle_{H_m(\mathbb B)}=\frac{(fu)_\alpha}{\rho_m(\alpha)}=\sum_{0\leq \beta\leq \alpha}\frac{f_\beta u_{\alpha -\beta}}{\rho_m(\alpha)}
\]
\[
=\langle u,\sum_{0\leq \beta \leq \alpha}\frac{\rho_m(\alpha - \beta)}{\rho_m(\alpha)}\overline{f}_\beta z^{\alpha - \beta}\rangle_{H_m(\mathbb B)}.
\]
Therefore the polynomials are contained in the domain of the adjoint $T^\ast_f$ of $T_f$ and
\[
T^\ast_f z^\alpha = \sum_{0\leq \beta \leq \alpha}\frac{\rho_m(\alpha - \beta)}{\rho_m(\alpha)}\overline{f}_\beta z^{\alpha - \beta}
\]
for $\alpha \in \mathbb N^n$.
In particular, for any fixed polynomial $p \in \mathbb C[z]$, the mapping
\[
H_m(\mathbb B) \rightarrow H_m(\mathbb B), f \mapsto T^*_f p
\]
is conjugate linear and continuous. 

We call a bounded linear operator $T \in L(H_m(\mathbb B))$ a Toeplitz operator with
pluriharmonic symbol $f$ if there are functions $g,h \in H_m(\mathbb B)$ with $f = g + \overline{h}$ and
\[
Tp = T_g p + T^*_h p  \quad \quad {\rm for \; all} \; p \in \mathbb C[z].
\]
An elementary argument shows that, although the representation $f = g + \overline{h}$  is only
unique up to an additive constant, the right-hand side of the above equation is independent of
the choice of $g$ and $h$. Furthermore, the function $f = g + \overline{h}$ is uniqely determined by $T$.
Indeed, if $g,h \in H_m(\mathbb B)$ satisfy the above relation with $T = 0$, then
\[
g + \overline{h(0)} = (T_g + T^*_h)(1) = 0.
\]
But then $T^*_h z^{\alpha} = - T_g z^{\alpha} = \overline{h(0)} z^{\alpha}$ for all $\alpha \in \mathbb N^n$.
Hence $h_{\beta} = 0$ for all $\beta \in \mathbb N^n \setminus \{ 0 \}$ and $g + \overline{h} = 0$. In the 
following we use the notation $T_f$ for the Toeplitz operator on $H_m(\mathbb B)$ with 
pluriharmonic symbol $f$.

\begin{thm} \label{bhimpliesph}
Let $T \in \mathcal T_{{\rm BH}}$ be given. Define
\[
g = (T - T_0)(1) \quad {\rm and} \quad h = T^*(1).
\]
Then $T = T_f$ is a Toeplitz operator with pluriharmonic symbol $f = g + \overline{h}$.
\end{thm}

\proof
For $k \in \mathbb Z$, let as before
\[
T_k = \frac{1}{2\pi}\int \limits^{2\pi}_0 e^{-ikt}U(t)TU(t)^* dt
\]
denote the $k$th homogeneous component of $T$. Define $g, h \in H_m(\mathbb B)$ and
the pluriharmonic function $f$ as in the statement of the theorem. Our aim is to
show that $T = T_f$. Set
\[
q_k = T_k1 \; \; {\rm for} \; \; k \geq 0, \quad q_k = \overline{(T_k)^*1} \; \; {\rm for} \; \; k < 0.
\]
Then
\[
T1 = \lim_{N \rightarrow \infty} \sum_{|k | \leq N} (1 -\frac{| k |}{N+1}) T_k1
   = \lim_{N \rightarrow \infty} \sum_{k=0}^N (1 -\frac{k}{N+1}) q_k
\]
and
\begin{align*}
T^*1 &= \lim_{N \rightarrow \infty} \sum^N_{k=0} (1 -\frac{k}{N+1}) (T^*)_k 1
   = \lim_{N \rightarrow \infty} \sum_{k=-N}^0 (1 -\frac{| k |}{N+1}) (T_k)^*1\\
	 &= \lim_{N \rightarrow \infty} \sum_{k=-N}^0 (1 -\frac{| k |}{N+1}) \overline{q}_k,
\end{align*}
where all sequences converge in $H_m(\mathbb B)$. Since $T, T^* \in \mathcal T_{{\rm BH}}$,
it follows from Lemma \ref{inherit} and Theorem \ref{multi} that $T_k = T_{T_k(1)} = T_{q_k}$
for $k \geq 0$ and that
\[
(T^*)_{-k} = T_{(T^*)_{-k}(1)} = T_{(T_k)^*1} = T_{\overline{q}_k}
\]
for $k < 0$. Let $p \in \mathbb C[z]$ be a polynomial. Because of $T_k = (T^*)_{-k}^*$ we find that
\[
Tp = \lim_{N \rightarrow \infty} \sum_{|k | \leq N} (1 -\frac{| k |}{N+1}) T_k p
\]
\[
= \lim_{N \rightarrow \infty} [ \sum_{k=0}^N (1 -\frac{k}{N+1}) q_k p \, - \, T_0(1)p \,
  + \sum_{k=-N}^0 (1 -\frac{| k |}{N+1}) T_{\overline{q}_k}^* p ]
\]
\[
= T_gp + \lim_{N \rightarrow \infty} \sum_{k=-N}^0 (1 -\frac{| k |}{N+1}) T_{\overline{q}_k}^* p.
\]
Since the mapping
$
H_m(\mathbb B) \rightarrow H_m(\mathbb B), u \mapsto T^*_up,
$
is conjugate linear and continuous, we conclude that
\[
Tp = T_g p + T^*_h p.
\]
Thus we have shown that $T = T_f$ with  $f = g + \overline{h}$.
\proofend

To prove that conversely each Toeplitz operator with pluriharmonic symbol satisfies
the Brown-Halmos condition, we use Lemma \ref{powers} and an approximation argument.

\begin{thm} \label{phimpliesbh}
Let $T = T_f \in L(H_m(\mathbb B))$ be a Toeplitz operator with pluriharmonic symbol
$f = g + \overline{h}$, where $g, h \in H_m(\mathbb B)$. Then $T \in \mathcal T_{{\rm BH}}$.
\end{thm}

\proof
Let $f = g + \overline{h}$ with $g, h \in H_m(\mathbb B)$. 
Let us denote by $g = \sum_{k=0}^{\infty} g_k$ and $h = \sum_{k=0}^{\infty} h_k$ the
homogeneous expansions of $g$ and $h$. Again using the continuity of the maps
\[
H_m(\mathbb B) \rightarrow H_m(\mathbb B), f \mapsto T^*_f p \quad (p \in \mathbb C[z])
\]
and the fact that $\delta M_z \mathbb C[z]^n \subset \mathbb C[z]$, we obtain as an
application of Lemma \ref{powers} that
\[
M'^*_z T_f M'_z (p_i) = \sum_{k=0}^{\infty} M^*_z \delta T_{g_k} \delta M_z (p_i)
+ \sum_{k=0}^{\infty} M^*_z \delta T^*_{h_k} \delta M_z (p_i)
\]
\[
= \sum_{k=0}^{\infty} P_{{\rm Im} M^*_z} \Big(\oplus \sum_{j=0}^{m-1} (-1)^j \binom{m}{j+1} 
\sigma^j_{M_z}(T_{g_k})\Big) P_{{\rm Im} M^*_z} (p_i)
\]
\[
+ \sum_{k=0}^{\infty} P_{{\rm Im} M^*_z} \Big(\oplus \sum_{j=0}^{m-1} (-1)^j \binom{m}{j+1} 
\sigma^j_{M_z}(T^*_{h_k})\Big) P_{{\rm Im} M^*_z} (p_i)
\]
for each tuple $(p_i) \in \mathbb C[z]^n$.
Using Lemma 1 in \cite{E} we find that
\[
( \oplus M_z^{* \alpha}) P_{{\rm Im} M^*_z} = (\oplus M_z^{* \alpha}) (M^*_z M_z)^{-1} (M^*_z M_z) = 
(\oplus M_z^{* \alpha}) M^*_z \delta M_z.
\]
This identity shows that the operator $(\oplus M_z^{* \alpha}) P_{{\rm Im} M^*_z}$ maps $\mathbb C[z]^n$ into itself.
Thus by reversing the above arguments we obtain that
\[
M'^*_z T_f M'_z (p_i) = P_{{\rm Im} M^*_z} \Big(\oplus \sum_{j=0}^{m-1} (-1)^j \binom{m}{j+1} 
\sigma^j_{M_z}(T_{g})\Big) P_{{\rm Im} M^*_z} (p_i)
\]
\[
+ P_{{\rm Im} M^*_z} \Big(\oplus \sum_{j=0}^{m-1} (-1)^j \binom{m}{j+1} 
\sigma^j_{M_z}(T^*_h)\Big) P_{{\rm Im} M^*_z} (p_i)
\]
\[
= P_{{\rm Im} M^*_z} \Big(\oplus \sum_{j=0}^{m-1} (-1)^j \binom{m}{j+1} \sigma^j_{M_z}(T_f)\big) P_{{\rm Im} M^*_z} (p_i)
\]
for each tuple $(p_i) \in \mathbb C[z]^n$. Since the polynomials are dense in $H_m(\mathbb B)$, the 
proof is complete.
\proofend

We finish by indicating that our definition of Toeplitz operators with pluriharmonic symbol coincides with
the usual one on the Hardy space and the weighted Bergman spaces.
Let $m \geq n$ be an integer and let $f:\mathbb B \rightarrow \mathbb C$ be a bounded pluriharmonic function. 
Then there are functions $g,h \in H_m(\mathbb B)$ with $f=g+\overline{h}$ (see for instance Proposition 6.1 in\cite{Zh}). 
Suppose first that $m \geq n+1$. Let $\mathcal T_f = P_{H_m(\mathbb B)} M_f|_{H_m(\mathbb B)}$, where
$P_{H_m(\mathbb B)}$ denotes the orthogonal projection of $L^2(\mathbb B,\mu_m)$ onto $H_m(\mathbb B)$
and $M_f$ is the operator of multiplication with $f$ on $L^2(\mathbb B,\mu_m)$ (cf. Section 1). Then
\[
\langle \mathcal T_f p,q \rangle = \langle fp,q \rangle = \langle gp,q\rangle + \langle p,hq \rangle
= \langle T_g p+T^\ast_h p,q \rangle
\]
for all polynomials $p,q \in \mathbb C [z]$. Next let us consider the case $m=n$. The boundary map
\[
h^\infty(\mathbb B)\rightarrow L^\infty(S),\ \varphi \mapsto \varphi^\ast
\]
defines an isometric isomorphism between the Banach space of all bounded $\mathcal M$-harmonic functions on $\mathbb B$ equipped with the supremum norm and $L^\infty(S)$ formed with respect to the normalized
surface measure on $S=\partial \mathbb B$. The inverse of this map is given by the Poisson transform (see Chapter 5 in \cite{Ru})
\[
L^\infty(S)\rightarrow h^\infty(\mathbb B),\ \varphi \mapsto P[\varphi].
\]
For $\varphi \in h^\infty(\mathbb B)$, the Toeplitz operator 
$\mathcal T_\varphi:H^2(\mathbb B)\rightarrow H^2(\mathbb B)$ is defined by
\[
\mathcal T_\varphi(u)=C[\varphi^\ast u^\ast],
\]
where the right-hand side denotes the Cauchy integral of $\varphi^* u^* \in L^2(S)$. 
For $f,g,h$ as above and any pair of polynomials $p,q \in \mathbb C[z]$, we obtain
\[
\langle \mathcal T_fp,q\rangle_{H^2(\mathbb B)} = 
\langle C[(gp)^\ast],q\rangle_{H^2(\mathbb B)}+\langle C[(\overline{h}p)^\ast],q\rangle_{H^2(\mathbb B)}.
\]
By Theorem 5.6.8 in \cite{Ru}, we have
\[
\langle C[(gp)^\ast],q\rangle_{H^2(\mathbb B)}=\langle gp,q\rangle_{H^2(\mathbb B)}
\]
and as an application of Theorem 5.6.9 in \cite{Ru} we obtain
\begin{align*}
&\langle C[(\overline{h}p)^\ast],q\rangle_{H^2(\mathbb B)}=\langle C[(\overline{h}p)^\ast,q^\ast\rangle_{L^2(S)}\\
&=\langle P_{H^2(S)}(\overline{h}p)^\ast,q^\ast\rangle_{L^2(S)}=\langle (\overline{h}p)^\ast, q^\ast\rangle_{L^2(S)}\\
&=\langle p^\ast,(hq)^\ast\rangle_{L^2(S)}=\langle p,hq\rangle_{H^2(\mathbb B)}.
\end{align*}
Thus for $m\geq n$, it follows that
\[
\langle\mathcal T_fp,q\rangle=\langle T_gp+T^\ast_hp,q\rangle
\]
for all polynomials $p,q \in \mathbb C[z]$. Hence $T_f = \mathcal T_f$ on $H_m(\mathbb B)$ for $m\geq n$.
In the particular case $n = 1 = m$ each Toeplitz operator $T_f \in L(H^2(\mathbb T))$ with
$L^{\infty}$- symbol $f \in L^{\infty}(S)$ coincides up to
unitary equivalence with a Toeplitz operator with harmonic symbol on $\mathbb D$, more precisely, 
$T_f \cong \mathcal T_{P[f]}$. Thus the results contained in Theorem 4 and Theorem 5 reduce to
the classical Brown-Halmos characterization \cite{BH} of Toeplitz operators on the Hardy space of the
unit disc in this case.

J. Eschmeier \\
Fachrichtung Mathematik\\
Universit\"at des Saarlandes\\
Postfach 151150\\
D-66041 Saarbr\"ucken\\
Germany\\[.7cm]

S. Langend\"orfer\\
Fachrichtung Mathematik\\
Universit\"at des Saarlandes\\
Postfach 151150\\
D-66041 Saarbr\"ucken\\
Germany

\texttt{eschmei@math.uni-sb.de, langendoerfer@math.uni-sb.de}

\end{document}